\documentclass[a4,amstex,12pt]{article}
\usepackage{amsmath,wrapfig}
\usepackage{amssymb,bm,epic}
\usepackage[mathscr]{eucal}
\usepackage{amsthm, amsfonts, latexsym,enumerate}
\usepackage{graphicx}
\usepackage{setspace}
\def\fixedfigure{\def\@captype{figure}}
\def\fixedtable{\def\@captype{table}}
\numberwithin{equation}{section}
\theoremstyle{definition}
\newtheorem{dfn}{Definition}[section]
\newtheorem{thm}[dfn]{Theorem}

\newtheorem{lmm}[dfn]{Lemma}

\newtheorem*{prf}{proof}
\newtheorem{alg}[dfn]{Algorithm}

\theoremstyle{break}
\begin{document}
\centerline{\Large Network Simplex Algorithm}
\vspace{2mm}
\centerline{\Large associated 
with the Maximum Flow Problem}\par\vspace{4mm}
\centerline{Sennosuke WATANABE$^a$, Hodaka TANAKA$^b$, 
Yoshihide WATANABE$^c$}\par\vspace{4mm}
\centerline{{\small $^a$ \textit{National Institute of Technology, 
Oyama College,}}} 
\par
\centerline{{\small \textit{771 Nakakuki, Oyama, Tochigi 323-0806 Japan}}}
\par
\centerline{{\small $^b$\textit{Graduate School of Science and 
Engineering,
Science of Environment and}}}
\par
\centerline{{\small \textit{Mathematical Modeling, 
Doshisha University, 1-3 Tatara 
Miyakodani,}}}
\par
\centerline{{\small \textit{Kyotanabe, 610-0394 Japan}}}
\par
\centerline{{\small $^c$\textit{Faculty of Science and Engineering, 
Department of Mathematical Sciences,}}}\par
\centerline{{\small \textit{Doshisha University, 1-3 Tatara Miyakodani, 
Kyotanabe, 610-0394 Japan}}}\par\vspace{2mm}
\centerline{{\small \textit{E-mail addresses:} 
$^a$sewatana@oyama-ct.ac.jp,}}\par 
\centerline{{\small $^c$yowatana@mail.doshisha.ac.jp}}
\begin{abstract}
In the present paper, 
we apply the network simplex algorithm 
for solving the minimum cost flow problem, 
to the maximum flow problem.
Then we prove that the cycling phenomenon
which causes the infinite loop in the algorithm,
does not occur in the network simplex 
algorithm associated with the maximum flow 
problem.
\end{abstract}
\textbf{Keywords:}
Network, Maximum Flow Problem, 
Network Simplex Algorithm, Cycling

\section{Introduction}
Two key optimization problems with respect to networks 
defined in directed graphs are the maximum flow problem 
(MFP) and the minimum cost flow problem (MCFP). 
These two problems have their own algorithms to 
achieve optimal solutions. 
It is well known that the MFP can be interpreted as 
a special class of the MCFP \cite{AMO,KV}.  
Thus, the network simplex algorithm for the minimum cost 
flow problem can be applied to the maximum flow problem. 
The network simplex algorithm is known to be one of the most 
efficient algorithms for solving the MCFP \cite{AMO,Jun}.  
However, a serious drawback of this algorithm is the so-called 
cycling phenomenon, which yields an infinite loop in 
the updating of feasible solutions, 
causing the algorithm to never reach an optimal solution. 
Thus, preventing this cycle from occurring is an indispensable part 
of effectively applying the algorithm. 
In the present paper, we will consider the MFP as a special class 
of the MCFP, and solve it by applying the network simplex algorithm. 
We will show that the network simplex algorithm for 
the MFP becomes very concise and simple. 
Furthermore, we will also demonstrate the main result of this paper, 
which is that the cycling phenomenon never occurs in this application.

\section{Minimum Cost Flow Problem and Maximum Flow Problem}
\subsection{Minimum Cost Flow Problem (MCFP)}
Let $G=(V,E)$ be a digraph with the vertex set $V$ and edge set $E$. 
As the upper and lower capacity of $e$, 
we will assign non-negative integers $b(e)$ and $c(e)$, respectively, 
to each edge $e\in E$ satisfying the inequality $0\leq b(e) \leq c(e)$. 
For the cost of $e$, we will assign a positive real number $\gamma(e)$ 
to each edge $e\in E$. 
Furthermore, for the demand of $v$, we will assign the integer 
$d(v)$ to each vertex $v$ with $\sum_{v\in V}d(v)=0$.
 
The quadruple $\mathcal{N}=(G,b,c,d)$ is called a network 
associated with the MCFP. 
A flow on the network $\mathcal{N}$ is the function $f$ on $E$ 
satisfying the following conditions (i) and (ii): 
\begin{enumerate}
\item[(i)] 
The capacity constraint at each edge: 
$$b(e) \leq f(e) \leq c(e) \quad 
\text{for all } e\in E \ . $$
\item[(ii)]
The demand condition at each vertex: 
$$
\sum_{\partial^+(e) = v} f(e) - \sum_{\partial^-(e) = v}
f(e) = d(v) \quad \text{for all } v \in V \ .
$$	 
\end{enumerate}
where the maps $\partial^-:E\to V$ and $\partial^+:E\to V$ with 
$\partial^-(e)=u$ and $\partial^+(e)=v$ for $e=(u,v)$. 
The MCFP is the problem of finding the flow $f$ whose cost 
$\gamma(f)=\sum_{e\in E} \gamma(e) f(e)$ is minimal. 

\subsection{Maximum Flow Problem (MFP)}
Let $G=(V,E)$ be a digraph with the vertex set $V$ and edge set $E$. 
The digraph $G$ has the source $s\in V$ and the sink $t\in V$, 
and each edge $e\in E$ has the capacity $c(e)$. 
Then, the quadruple $\mathcal{N}=(G,c,s,t)$ is called a network 
associated with the MFP. 
A flow on the network $\mathcal{N}$ is the function $f$ on $E$ 
satisfying the following conditions (i) and (ii):
\begin{enumerate}
\item[(i)] The capacity constraint at each edge:
$$0\leq f(e) \leq c(e) \quad \text{for all } e\in E \ .$$ 
\item[(ii)] The flow conservation law 
at each vertex except for the source $s$ and the sink $t$:
$$
\sum_{\partial^+(e) = v} f(e) - \sum_{\partial^-(e) = v}f(e) 
=0 
\quad \text{for all } v \in V \backslash \{s,t\} \ .
$$
\end{enumerate}
It follows from the flow conservation law that 
$$
\tau(f)=\sum_{\partial^-(e)=s}f(e) = \sum_{\partial^+(e)=t}f(e) \ .
$$
where $\tau(f)$ is the flow value of the flow $f$. 
The MFP is the problem of finding the flow such that 
the flow value is the maximum.

\subsection{Minimum Cost Flow Problem and Maximum Flow Problem}
In this subsection, we will explain how the MFP can be interpreted 
as a special case of the MCFP. 
Let $\mathcal{N}=(G,c,s,t)$ be a network on the digraph $G=(V,E)$ 
associated with the MFP, 
and let $\tau(f)$ be the flow value of a flow $f$ on $\mathcal{N}$. 
By adding a new edge $r =(t,s)$ to $E$, 
giving $E'=E\cup \{r\}$, we can define 
the flow network $\mathcal{N}'=(G',b',c',d')$ on 
the digraph $G'=(V,E')$ associated with the MCFP as follows: 
\begin{enumerate}
\item[(i)] 
$b'(e)=0$ for all $e\in E'$ and $c'(e)=c(e)$ for $e\in E$ 
and $c'(r)=M (\text{sufficiently large})$.
\item[(ii)]
$d'(v)=0$ for all $v\in V$. 
\end{enumerate}
Now, let $f$ be a function on the edge set $E$. 
We will extend $f$ to the function $f'$ on the edge set $E'$ 
using $f'(r)=\tau(f)$. 
Then, $f$ is a flow on the network $\mathcal{N}=(G,c,s,t)$ 
if and only if $f'$ is a flow on the network 
$\mathcal{N}'=(G',b',c',d')$. 
Furthermore, we will define the cost $\gamma$ on $E'$ 
using $\gamma(e)=-1$ for $\partial^+(e)=t$ and $\gamma(e)=0$ otherwise. 
Then, a flow $f$ on $\mathcal{N}$ is the maximum flow 
if and only if the flow $f'$ is the minimum cost flow 
on $\mathcal{N}'$ with respect to the cost $\gamma$.

\section{Network Simplex Algorithm}
The network simplex algorithm consists of the following 3 steps: 
\begin{enumerate}
\item
Find the initial feasible tree structure.
\item 
Decide whether or not the given tree structure is optimal.
\item 
If the tree structure is not optimal, 
then update the tree structure and associated tress solution 
to improve the cost.
\end{enumerate}
In the present paper, we will refer the reader to \cite{Jun} 
for a description and use of step 1, which we do not address.

\subsection{Tree Structure}
We will consider the MCFP on a digraph $G=(V,E)$. 
Let $T=(V,E_T)$ be a spanning tree of $G$. 
We often identify the spanning tree $T$ and its edge set $E_T$ 
and use the same symbol $T$ for $E_T$.  
We will divide the edge set $E \setminus T$ into 
the disjoint union $E\setminus T = L \cup U$, 
where subsets $L,U\subset E$ are allowed to be empty. 
We can fix the triplet $(T,L,U)$ and define the function $f$ 
on the edge set $E$ as follows: 
First, we define the function $f$ on $L$ and $U$ using 
$f(e)=c(e)\;(e\in U)$ and $f(e)=b(e)\;(e\in L)$, respectively. 
Then, the demand condition for the function $f$ uniquely determines 
the flow value on the spanning tree $T$. 
Thus, we can see that the function $f$ on $E$ is uniquely determined 
by the triplet $(T,L,U)$. 
Note that the function $f$ does not always become 
the flow on the network $\mathcal{N}$ because it does not 
always satisfy the capacity constraints on $T$. 
If the function $f$ uniquely determined by the triplet $(T,L,U)$ 
becomes the flow on the network, the flow $f$ is called 
the tree solution associated with the feasible 
tree structure $(T,L,U)$. 
We call an edge $e$ free with respect to the flow on 
the network $\mathcal{N}$ provided that $b(e)<f(e)<c(e)$ holds. 
It follows directly from the definition of the tree solution 
that the free edge with respect to the tree solution must be 
contained in $T$. 
However, all edges of $T$ need not be free with respect to 
the tree solution. 
A feasible tree structure is called optimal if 
the unique tree solution associated with the tree structure is optimal. 
It is well known that if the MCFP has an optimal solution, 
then there exists an optimal tree solution \cite{AMO,Jun}.  
The network simplex algorithm is an algorithm for solving 
the MCFP by updating the feasible tree structure.

\subsection{Optimality Condition}
Consider the MCFP on the network $\mathcal{N}$ defined on 
a digraph $G$ and let $f$ be the tree solution associated with 
a feasible tree structure $(T,L,U)$. 
We can determine the potential $\pi:V \to \mathbb{R}$ using 
the following procedure: 
Fix a base vertex $x \in V$, and set $\pi(x)=0$. 
For each vertex $v \in V \setminus \{x\}$, 
there exists a unique undirected path $P_v$ from $x$ to $v$ 
in the spanning tree $T$. 
We will define the direction of the path $P_v$ from $x$ to $v$. 
Then, all the edges in $P_v$ are categorized as forward and 
backward edges with respect to the direction of $P_v$; 
the set of forward edges is denoted by $P_v^+$ and 
the set of backward edges by $P_v^-$. 
We can define the potential $\pi$ using 
$$
\pi(v)= \sum_{e\in P_{v}^+} \gamma(e) 
- \sum_{e \in P_v^-} \gamma(e)
\quad \text{for all } v \in V \setminus \{x\} \ .
$$
Furthermore, we can define 
the reduced cost $\gamma_{\pi}:E \to \mathbb{R}$ 
in terms of the potential $\pi$ using 
$$
\gamma_{\pi}(e) = \gamma(e)+\pi(v_i)-\pi(v_j) 
\quad \text{for all } e=(v_i,v_j)\in E \ .
$$
Then, we have 
$$
\gamma_{\pi}(f) = \sum_{e\in E} f(e)\gamma_{\pi}(e) 
= \gamma(f) - \sum_{v\in V}\pi(v)d(v) \ ,
$$
which shows that the cost $\gamma(f)$ and 
the reduced cost $\gamma_{\pi}(f)$ of the tree solution $f$ 
differ only by a constant independent. 
Thus, it is easy to see that a tree structure $(T,L,U)$ is optimal 
if and only if the reduced cost $\gamma_{\pi}$ on $(T,L,U)$ 
satisfies the following condition:
\begin{align}\label{Opti}
\gamma_{\pi} (e)
\left\{
\begin{array}{lll}
= 0    & \text{ for all } & e\in T \ ,\\
\geq 0 & \text{ for all } & e\in L \ ,\\
\leq 0 & \text{ for all } & e\in U \ .
\end{array}
\right.
\end{align}

\subsection{Network Simplex Algorithm}
\begin{alg}\label{NSA}
Consider the MCFP on a digraph $G$.
We can obtain the optimal solution of the MCFP 
using the following procedures:
\begin{enumerate}
\item[1.] 
If the reduced cost $\gamma_{\pi}$ with respect to the feasible tree
structure $(T,L,U)$ 
satisfies optimality condition (\ref{Opti}),  
then the tree solution associated with the tree structure $(T,L,U)$ is 
the optimal solution of the MCFP.
Otherwise go to the next step.
\item[2.]
Choose some $e \in L \cup U$ 
satisfying either 
\begin{enumerate}
\item[(i)] 
$e \in L \text{ and } \gamma_{\pi}(e)< 0$ or 
\item[(ii)] 
$e \in U \text{ and } \gamma_{\pi}(e)> 0$.  
\end{enumerate}
We call this edge $e$ the entering edge.
Then, the edge set $T\cup \{e \}$ contains 
the unique circuit $C_T(e)$. 
We can determine the orientation of 
$C_T(e)$ as the direction of the entering edge $e$ if $e\in L$, 
and as the opposite direction of $e$ if $e\in U$.
\item[3.]
Augment the flow $f$ along the circuit $C_T(e)$ until either 
\begin{enumerate}
\item[(i)] 
one of the forward edges in $C_T(e)$ reaches 
its upper capacity or 
\item[(ii)] 
one of the backward edge in $C_T(e)$ reaches 
its lower capacity. 
\end{enumerate}
Choose an edge $a$ in $C_T(e)$ that reaches the upper or 
lower capacity and call $a$ the leaving edge.   
\item[4.]
Update the new tree structure $(T',L',U')$ as follows:
\begin{align*}
T'&:= (T\cup \{e\})\setminus \{a \} , \\
L'&:= 
\left\{
\begin{array}{ll}
(L\setminus \{e \})\cup \{a \} & 
\text{ if } a \text{ reaches its lower capacity bound} \ ,\\
L\setminus \{e \} & 
\text{ if } a \text{ reaches its upper capacity bound} ,
\end{array}
\right. \\
U'&:= E\setminus (T\cup L) \ .
\end{align*}
Compute the potential $\pi$ of $(T',L',U')$ and 
the reduced cost $\gamma_{\pi}$.
Go to Step 1.
\end{enumerate}
\end{alg}

\section{Network Simplex Algorithm associated with the Maximum Flow Problem}
In this section, 
we will apply the network simplex algorithm to the MFP.

\subsection{Pseudo Tree Structure}
We will consider the MFP on a network $\mathcal{N}=(G,c,s,t)$ 
with the digraph $G=(V,E)$. 
Let $A=(V_A,E_A)$ and $B=(V_B,E_B)$ be trees of $G$ satisfying 
$$
s\in V_A , t\in V_B, V_A\cup V_B =V \text{ and } 
V_A\cap V_B =\emptyset \ .
$$
We will express the union of $E_A$ and $E_B$ as $T$. 
Furthermore, we will divide the edge set $E \setminus T$ 
into the disjoint union $E\setminus T = L \cup U$. 
We can fix the triplet $(T,L,U)$ and define 
the function $f$ on the edge set $E$ as follows: 
Define the function $f$ on $L$ and $U$ using 
$f(e)=0\;(e\in L)$ and $f(e)=c(e)\;(e\in U)$. 
Then, the values of the function $f$ on $T$ are uniquely 
determined by the conservation law. 
Thus, we see that the function $f$ on $E$ is uniquely 
determined by the triplet $(T,L,U)$. 
If this function $f$ becomes the flow on the network, 
it is called the tree solution associated with 
the feasible pseudo tree structure $(T,L,U)$. 

\subsection{Initial Feasible Pseudo Tree Structure}
Let $\mathcal{N}=(G,c,s,t)$. 
Consider a pseudo tree structure $(T,L,U)$ 
with $L=E\setminus T$ and $U=\emptyset$.
Then the function $f$ which is uniquely determined by $(T,L,U)$ 
becomes the zero flow on the network $\mathcal{N}$.
We can take the zero flow $f$ as the initial tree solution associated 
with the initial feasible pseudo tree structure $(T,L,U)$. 

\subsection{Optimality Condition}
Consider the MFP on the network $\mathcal{N}$ 
with the digraph $G$ and 
let $f$ be the tree solution associated with a feasible 
pseudo tree structure $(T,L,U)$.
We will define sets of edges $L_{AB}$ and $U_{BA}$ as 
$$
L_{AB}=\{(u,v)\in L|u\in V_A, v\in V_B \} \text{ and } 
U_{BA}=\{(v,u)\in U|u\in V_A, v\in V_B \} \ .
$$
Then the optimality condition of the tree structure 
is given as follows: 
\begin{lmm}\label{lem}
The tree solution $f$ associated with a feasible 
pseudo tree structure $(T,L,U)$ is optimal if and only if 
sets of edges $L_{AB}$ and $U_{BA}$ are empty sets, that is, 
$L_{AB}=\emptyset$ and $U_{BA}=\emptyset$.
\end{lmm}
\begin{prf}
The tree solution associated with 
a pseudo tree structure $(T,L,U)$ on the network 
associated with the MFP accords with 
a tree solution associated with the tree structure 
$(T\cup \{(t,s)\},L,U)$ on the network associated with 
MCFP. 
Thus, we only have to prove that the edge set given in (\ref{Opti})
that obstruct the optimality of the solution coincides with the edge 
set $L_{AB}\cup U_{BA}$ of Lemma \ref{lem}.
We will choose an arbitrary feasible pseudo tree structure $(T,L,U)$ 
on the network associated with the MFP. 
We will divide the edge set $E\setminus T$ into six subsets:
\begin{align*}
L_{AB}&= \{(u,v)\in L\ |\ u\in V_A,v\in V_B\}, &
U_{BA}&=\{(v,u)\in U\ |\ u\in V_A,v\in V_B\}\\
E_{AA}&=\{(u,v)\in E\setminus T\ |\ u,v\in V_A\},   & 
E_{BB}&=\{(u,v)\in E\setminus T\ |\ u,v\in V_B\} \\
L_{BA}&=\{(v,u)\in L\ |\ u\in V_A,v\in V_B\}, & 
U_{AB}&=\{(u,v)\in U\ |\ u\in V_A,v\in V_B\}
\end{align*}
First, we will prove that the edges in $L_{AB}$ or $U_{BA}$ do not 
satisfy optimality condition \eqref{Opti}. 
To compute the potential $\pi$, we will take the sink vertex $t$ as the 
base vertex and set $\pi(t)=0$. 
\\
{\bf (1-i)} $(u,v)\in L_{AB}$ 
\\
For $u\in V_A$ we have $\pi(u)=0$ and for $v\in V_B$, we have
$\pi(t)=0$ and $\pi(v)=1$ for $v\neq t$. 
Since $\gamma(u,t)=-1$ and $\gamma(u,v)=0$ for  $v\not= t$,  
we can compute $\gamma_\pi(u,t)=-1$ and 
$\gamma_\pi(u,v) =  -1$ for $v\not= t$. 
Thus, $\gamma_\pi(u,v)=-1$ for all $(u,v)\in L_{AB}$, 
which shows that 
$(u,v)\in L_{AB}$ does not satisfy optimality condition \eqref{Opti} 
(cf. Figure \ref{1-i,ii}).
\\
{\bf (1-ii)} $(v,u)\in U_{BA}$
\\
It follows from the assumption on the MFP that 
there are no edges with tail $t$. 
We have $\pi(u)=0$ for  $u\in V_A$ and
$\pi(v)=1$ for $v\in V_B$.  
Since $\gamma(v,u)=0$ for all $(v,u)\in U_{BA}$, 
we can see that $\gamma_{\pi}(v,u)=0+1-0=1$ 
(cf. Figure \ref{1-i,ii}), which shows that $(v,u)\in U_{BA}$ 
does not satisfy optimality condition (\ref{Opti}).
\newpage
\begin{figure}[htbp]
\begin{center}
\includegraphics[width=135mm]{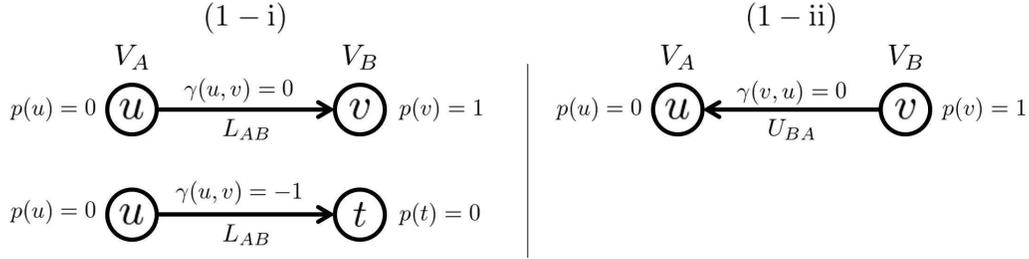}
\end{center}
\caption{The potential and the cost in the case of (1-i),(1-ii)} 
\label{1-i,ii} 
\end{figure}
Next, we will prove that edges in $E_{AA},E_{BB},L_{BA}$ and $U_{AB}$ satisfy 
optimality condition (\ref{Opti}).
\\
{\bf (2-i)} $(u,v)\in E_{AA}$
\\
Since $\pi(u)=\pi(v)=0$ and $\gamma(u,v)=0$ for $u,v \in V_A$, 
we can compute $\gamma_{\pi}(u,v)=0$ (cf. Figure \ref{2-i,ii}), 
which shows that $(u,v)\in E_{AA}$ satisfies optimality condition 
(\ref{Opti}).
\\
{\bf(2-ii)} $(u,v)\in E_{BB}$
\\
If $u\not= t,v\not= t$, 
we have $\pi(u)=\pi(v)=1$ and $\gamma (u,v)=0$ for  $u,v\in V_B$.
Therefore, we compute $\gamma_{\pi}(u,v)=0+1-1=0$ (cf. Figure \ref{2-i,ii}), 
which shows that $(u,v) \in E_{BB}$ satisfies optimality condition
(\ref{Opti}).
If $u\not= t,v=t$, we have $\pi(t)=0,\pi(u)=1$ and
$\gamma (u,v)=-1$.
Therefore, we compute $\gamma_{\pi}(u,v)=-1+1+0=0$ (cf. Figure \ref{2-i,ii}), 
which shows that $(u,v) \in E_{BB}$ satisfies the optimality condition
(\ref{Opti}).
\begin{figure}[htbp]
\begin{center}
\includegraphics[width=135mm]{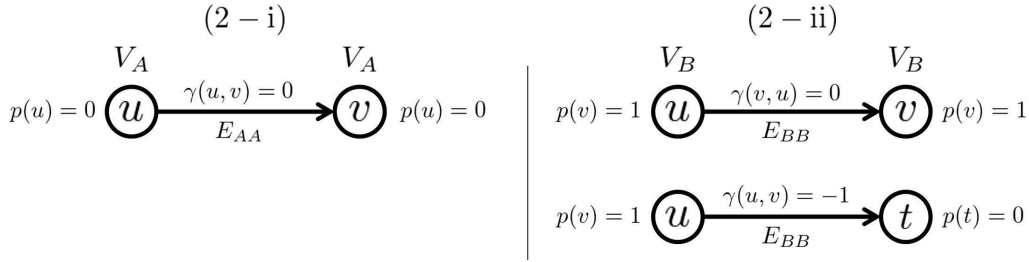}
\end{center}
\caption{The potential and the cost in the case of (2-i),(2-ii)}
\label{2-i,ii}
\end{figure}
\\
{\bf(2-iii)} $(v,u)\in L_{BA}$
\\
It follows from the assumption on the MFP that 
there are no edges with head $s$ or edges with tail $t$. 
Thus, $\pi(u)=0$ for $u\in V_A$ and $\pi(v)=1$ for $v\in V_B$. 
Furthermore, $\gamma (u,v)=0$ for $(v,u)\in L_{BA}$. 
Therefore, we can compute $\gamma_{\pi}(u,v)=0+1-0=1$ (cf. Figure \ref{2-iii,iv}), 
which shows that $(v,u) \in L_{BA}$ satisfies optimality condition 
(\ref{Opti}).
\\
{\bf(2-iv)} $(u,v)\in U_{AB}$
\\
We have $\pi(u)=0$ for $u\in V_A$.  
Furthermore, for $v\not= t$, we have $\pi(t)=0$ and $\pi(v)=1$, 
and $\gamma(u,t)=-1$ and $\gamma(u,v)=0$.  
Therefore, we can compute the reduced
cost $\gamma_\pi(u,v)=-1$ for both cases $v=t$ and  $v\neq t$  
(cf. Figure \ref{2-iii,iv}), which shows that $(u,v)\in U_{AB}$ 
satisfies optimality condition (\ref{Opti}).
\\
\begin{figure}[htbp]
\begin{center}
\includegraphics[width=135mm]{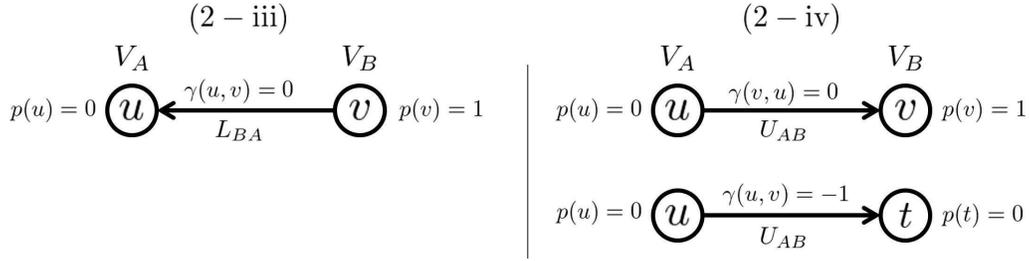}
\end{center}
\caption{The potential and the cost in the case of (2-iii),(2-iv)}
\label{2-iii,iv}
\end{figure}
\\
Thus, we have completed the proof of Lemma \ref{lem}. 
\qed
\end{prf}

\subsection{Network Simplex Algorithm associated 
the Maximum Flow Problem}
Consider the MFP on the network $\mathcal{N}=(G,c,s,t)$ 
on the digraph $G=(V,A)$.
We can describe the network simplex algorithm for MFP as follows.
\begin{alg} \ \\
{\bf Step 0:}
Choose arbitrary trees $A=(V_A,E_A)$ and $B=(V_B,E_B)$ 
satisfying 
$$
s\in V_A , t\in V_B , V_A \cup V_B \text{ and } 
V_A \cap V_B \ .
$$
Define the initial pseudo tree structure $(T,L,U)$ 
as $T=E_A\cup E_B$, $L=E\setminus T$ and $U=\emptyset$.
Take the zero flow as the initial tree solution 
associated with the tree structure $(T,L,U)$. 
\\
{\bf Step 1:}
If the pseudo tree structure $(T,L,U)$ satisfies 
$L_{AB} \cup U_{BA}=\emptyset$, 
then the tree solution $f$ associated with $(T,L,U)$ is optimal.
Else go to step 2.
\\
{\bf Step 2:}
Chose an edge $e\in L_{AB}\cup U_{BA}$ and attach $e$ to  $T$.
Then, there is the unique $s$-$t$ path $P$ in $T$.  
Determine the orientation of the path $P$ as the direction of $e$ if 
$e\in L_{AB}$ and as the opposite direction of $e$ if $e\in U_{BA}$.
Augment $f$ until at least one of the forward edges of $P$ reaches its upper 
capacity or one of the backward edge reaches its lower capacity, 
and choose one edge $a$ that is saturated 
by the above augmentation of the flow. 
\\
{\bf Step 3:}
Update the new pseudo tree structure $(T',L',U')$ as follows:
\begin{align*}
T'&:= (T\cup \{e\})\setminus \{a\} , \\
L'&:= 
\left\{
\begin{array}{ll}
(L\setminus \{e \})\cup \{a \} & 
\text{ if } a \text{ reaches its lower capacity bound} \ ,\\
L\setminus \{e \} & 
\text{ if } a \text{ reaches its upper capacity bound} ,
\end{array}
\right. \\
U'&:= E\setminus (T\cup L) \ .
\end{align*}
Go to Step 1.
\end{alg}

\section{Main Result - no cycling phenomenon}
The network simplex algorithm \ref{NSA} for the MCFP 
does not necessarily terminate in a finite number of iterations.
This inconvenience is caused by the degeneracy of 
the tree structure.
If the spanning tree $T$ in the tree structure $(T,L,U)$  
contains non-free edges, then the tree structure $(T,L,U)$ 
is called degenerate.
In this case, a proper augmentation along the 
circuit determined in Step 2 of algorithm \ref{NSA} 
may be impossible. In such a case only the feasible tree structure 
is updated and the tree solution associated with the tree structure 
is not updated. 
In this case, we may return to the same tree structure 
after a number of iterations without updating the tree solution.   
This phenomenon is known as cycling.
It is known that this cycling can be prevented by choosing 
the leaving edge appropriately, for example, 
by applying the so-called 
rule of the last blocking edge or the rule of the first blocking edge 
\cite{Jun,Senn1}.
In the present paper, we will prove that 
the cycling phenomenon never occurs in 
the network simplex algorithm for the maximum flow problem.
\begin{thm}
The cycling phenomenon never occurs in 
the network simplex algorithm for the maximum flow problem.
\end{thm}
\begin{prf}
We assume on the contrary that the cycling occurs and prove
that this leads to a contradiction. 
We will start with some pseudo tree structure $(T,L,U)$, and
return to the original pseudo structure $(T,L,U)$ 
without updating the flow value 
after a finite number of iterations. 
Let $A=(V_A,E_A)$ and $B=(V_B,E_B)$ be the corresponding trees 
associated with the pseudo tree structure $(T,L,U)$.  
We assume that there exists an 
edge $e_0 \in E$ that leaves the tree $T$ and later re-enters the tree. 
We call such an edge a leaving-entering edge. 
We will prove that the existence of 
leaving-entering edge gives rise to a contradiction. 
\par
First, we will consider 
the case where the leaving-entering edge $e_0$ is in $E_A$. 
If we consider the tree $A$ as an undirected graph,
there exists a unique path in $A$ between the source vertex $s$ and 
an arbitrary vertex $v\in V_A$ such that the vertices in $V_A$ can be layered 
as $V_A = \cup_{i\geq 0}V_i$, where $V_0=\{s\}$ and 
$V_i$ is the set of vertices such that the length of the unique path 
from $s$ in $A$ is $i$. 
Also, we define the set of edges $ E_j\subset E_A$ 
using $E_j =E_j^{f}\cup E_j^b$, where we set 
$$
E_j^f=\{(u,v)\in E_A| u\in V_{j-1},v\in V_j\},\quad
E_j^b=\{(v,u)\in E_A|u\in V_{j-1},v\in V_j\} \ .
$$ 
We assume $e_0\in E_k$ and that $k$ is taken as a minimum. 
Since the edge $e_0=(u,v)$ leaves from $E_A$, 
we can see that 
$f(e_0)=c(e_0)$ if $e_0\in E_k^f$ or $f(e_0)=0$ if  
$e_0\in E_k^b$. 
It follows from the minimality of $k$ that  
$V_0,V_1,\dots,V_{k-1}\subset V_A$ and 
$E_0,E_1,\dots,E_{k-1}\subset E_A$. 
For the edge $e_0$ to re-enter $E_A$, 
we must have either $e_0=(u,v)\in U_{BA}$ if $e_0\in E_k^f$ 
or $e_0=(v,u)\in L_{AB}$ if $e_0\in E_k^b$. 
In both cases, 
we have $u\in V_B$, which contradicts 
the assumption $u\in V_{k-1}\subset V_A$. 
We can also consider the assumption that 
there exists a leaving-entering edge $e_0\in E_A$, 
which similarly leads to a contradiction.
Thus, the proof of the main theorem is complete.
\qed
\end{prf}



\begin{thebibliography}{Abc}
\bibitem[1]{AMO} K.~Ahuja, L.~Magnanti and B.~Orlin:
\textit{NETWORK FLOWS},
Prentice-Hall, New Jersey, 1993.
\bibitem[2]{BK} A.~Bachem and W.~Kern: 
\textit{Linear Programming Duality}, 
Springer Verlag, Berlin Heidelberg, 1992.
\bibitem[3]{Jun}D.~Jungnickel: 
\textit{Graphs, Networks and Algorithm}, 
Second Edition, Springer Verlag, Berlin Heidelberg, 2005.
\bibitem[4]{KV}B.~Korte and J.~Vygen: 
\textit{Combinatorial Optimization},
Second Edition, Springer Verlag, 2001.
\bibitem[5]{Senn1}S.~Yoneda, S.~Watanabe and Y.~Watanabe: 
\textit{A New Rule for Preventing the Cycling in the Network 
Simplex Algorithm}, 
The Science and Engineering Review of Doshisha University, 
Science and Engineering Research Institute of Doshisha University, 
Vol.53 No.1, pp.54-57, 2012 (in Japanese).
\end{thebibliography}
\end{document}